\documentclass[11pt]{article}
\usepackage{amsmath}
\usepackage{amssymb}
\usepackage{amsthm}
\usepackage[stdtext]{youngtab}
\newcommand\ba{\bigskipamount}
\newcommand\bLP{\\[\ba]}
\newcommand\bPP{\\[\ba]\indent}

\newcommand\CC{\mathbb{C}}
\newcommand\RR{\mathbb{R}}
\newcommand\TT{\mathbb{T}}
\newcommand\ZZ{\mathbb{Z}}
\newcommand\al\alpha
\newcommand\be\beta
\newcommand\de\delta
\newcommand\tha\theta
\newcommand\la\lambda
\newcommand\De{\Delta}

\newcommand\half{\frac12}
\newcommand\thalf{\tfrac12}
\newcommand\iy\infty
\newcommand\wt{\widetilde}
\newcommand\lan{\langle}
\newcommand\ran{\rangle}
\newcommand{\hyp}[5]{\,\mbox{}_{#1}F_{#2}\!\left(
  \genfrac{}{}{0pt}{}{#3}{#4};#5\right)}
\newcommand{\qhyp}[5]{\,\mbox{}_{#1}\phi_{#2}\!\left(
  \genfrac{}{}{0pt}{}{#3}{#4};#5\right)}
 \newcommand\LHS{left-hand side}
\newcommand\RHS{right-hand side}

\numberwithin{equation}{section}
\newtheorem{theorem}{Theorem}[section]

\newtheorem{Definition}[theorem]{Definition}
\newenvironment{definition}{\begin{Definition}\rm}{\end{Definition}}
\newtheorem{Remark}[theorem]{Remark}
\newenvironment{remark}{\begin{Remark}\rm}{\end{Remark}}
\newtheorem{Example}[theorem]{Example}

\newcommand\Proof{\noindent{\bf Proof}\quad}

\begin{document}
\title{Okounkov's $BC$-type interpolation Macdonald polynomials
and their $q=1$ limit}
\author{Tom H. Koornwinder,
{\small {\tt T.H.Koornwinder@uva.nl}}}
\date{}
\maketitle
\begin{abstract}
This paper surveys eight classes of polynomials associated with $A$-type
and $BC$-type root systems: Jack, Jacobi, Macdonald and Koornwinder
polynomials and interpolation (or shifted) Jack and Macdonald polynomials
and their $BC$-type
extensions. Among these the $BC$-type interpolation Jack polynomials were
probably unobserved until now. Much emphasis is put on combinatorial formulas
and
binomial formulas for (most of) these polynomials. Possibly new results derived
from these formulas are a limit from Koornwinder to
\mbox{Macdonald} polynomials,
an explicit formula for Koornwinder polynomials in two variables, and
a combinatorial expression for the coefficients of the expansion of $BC$-type
Jacobi polynomials in terms of Jack polynomials which is different
from Macdonald's combinatorial expression. For these last coefficients
in the two-variable case the explicit expression
in Koornwinder \& Sprinkhuizen (1978)
is now obtained in a quite different way.
\end{abstract}
%
%
\section{Introduction}
In the past half century special functions associated with root
systems became an active area of research with many interconnections
and applications.  The early results were strongly motivated by the
notion of spherical functions on Riemannian symmetric spaces. An
ambitious program, which still has not come to an end, started to do
``zonal spherical analysis'' without underlying group and for a wider
parameter range than the discrete set of parameter values for which a
group theoretic interpretation is possible.  Another motivation came
from applications in multivariate statistics.
By the end of the eighties of the past century Heckman and Opdam
consolidated the
theory of Jacobi polynomials associated with root systems. In the same
period Macdonald,
in his {\em annus mirabilis} 1987, introduced the $q$-analogues of
these Jacobi polynomials in several manuscripts which were circulated
and eventually published: {\em Macdonald polynomials} $P_\la(x;q,t)$
(associated with $A$-type root systems) in \cite{28} and
\cite[Ch.~VI]{5},
{\em Macdonald polynomials associated with root systems} in \cite{21},
and scratch notes about hypergeometric
functions (associated with $BC$-type root systems) in \cite{3}. Again
in the same period Dunkl introduced his {\em Dunkl operators}, which
inspired Heckman, Opdam and in particular Cherednik to consider the
Weyl group invariant ($W$-invariant) special functions as part of a
more general theory of non-symmetric special functions which are
eigenfunctions of operators having a reflection term. Special
representations of graded and double affine Hecke algebras (DAHA's)
were an important tool. This approach not only introduced new
interesting special functions, but also greatly simplified the
$W$-invariant theory.

The author \cite{6} introduced a 5-parameter class of $q$-polynomials,
on the one hand extending the
3-parameter class of Macdonald polynomials associated with root system $BC_n$
\cite{21} and on the other hand providing the $n$-variable analogue of the
Askey-Wilson polynomials \cite{31}.
These polynomials became known in literature
as {\em Macdonald-Koornwinder} or
{\em Koornwinder polynomials}.
Cherednik's DAHA approach could also be used for these
polynomials, see Sahi \cite{18}, \cite{38} and Macdonald's monograph \cite{20}.
A different approach started by work of Sahi, Knop, Okounkov and Olshanski
(\cite{37}, \cite{15}, \cite{13}, \cite{14}, \cite{19}, \cite{11}, \cite{1}).
It used the so-called
{\em shifted} or {\em interpolation} versions
of Jack and Macdonald polynomials.
These could be characterized very briefly by their vanishing property
at a finite part of a \mbox{($q$-)}lattice, they could be represented
by {\em combinatorial formulas} (tableau sums) generalizing those for
Jack and Macdonald polynomials, and they occurred in generalized
{\em binomial formulas}. In particular, Okounkov's \cite{1}
{\em $BC_n$ type interpolation Macdonald polynomials} inspired
Rains \cite{2} to
use these in the definition of Koornwinder polynomials, thus building
the theory of these latter polynomials in a completely new way. An
analogous approach then enabled Rains to develop a theory of
{\em elliptic} analogues of Koornwinder polynomials, as surveyed in
\cite{30}.

Jack and Macdonald polynomials in $n$ variables play a double role, on
the one hand as homogeneous orthogonal polynomials associated with
root system $A_{n-1}$, on the other hand as generalized ``monomials''
(in the one-variable case ordinary monomials) in terms of which
orthogonal polynomials associated with root system $BC_n$ can be
naturally expanded.
This second role is emphasized in the approach using interpolation
polynomials, in particular where it concerns binomial formulas.

The present paper surveys, mainly in Sections 4 and 5
and after some preliminaries in Section~3,
the definition and properties of eight classes of
polynomials:
four associated with root system $BC_n$ and four with root system
$A_{n-1}$. Also four of these classes are for general $q$ and four are
for $q=1$.
Four of these classes can be considered as orthogonal polynomials
while the other four (interpolation) classes only play a role as
generalized monomials.  There are many limit connections between these
eight classes. For six of them (however, see \cite{36} and Remark
\ref{103}) combinatorial formulas are known, see such formulas mainly
in Section 6.
In a sense these combinatorial formulas are generalized hypergeometric series. 

One of the eight classes, the $BC_n$-type interpolation Jack polynomials,
seems to have been overlooked in literature, although it occurs very
naturally in the scheme formed by the limit connections.
It will be defined in Section 7. All its properties will be obtained
here as limit cases of properties of $BC_n$-type interpolation
Macdonald polynomials, including the combinatorial formula for
polynomials of this latter class.

Binomial formulas as they were already known
for three classes of polynomials
are surveyed in Section 8. The probably new binomial formula for
$BC_n$-type interpolation Jack polynomials is given in Section 9.
It gives a new approach to coefficients of the expansion of
$BC_n$-type Jacobi polynomials in terms of Jack polynomials.
As a consequence of the binomial formulas
a new limit formula \eqref{56} and a new proof of
an already known limit formula \eqref{62} will follow.

All classes of polynomials and formulas for them become much more
elementary and explicit in the one-variable case. This is the subject
of the Prelude in Section 2.  The two-variable case is already more
challenging, but explicit formulas are feasible.  This is the topic of
Sections 10 and 11. In particular, in Subsection 11.2 we arrive at an
explicit expression for $BC_2$-type Jacobi polynomials which was
earlier obtained in a very different way by the author together with
Sprinkhuizen in \cite{4}.
\bLP
{\em Acknowledgement}\quad
The material of this paper was first presented in lectures given at the
72nd S\'eminaire Lotharingien de Combinatoire,
Lyon, France, 24--26 March 2014. I thank the organizers for the invitation.
I thank Siddhartha Sahi, Ole Warnaar, Genkai Zhang and an anonymous
referee for helpful remarks. Thanks also to Masatoshi Noumi for making
available to me his unpublished slides on interpolation functions
of type $BC$.
\bLP
{\em Notation}\quad See \cite{29}.
Throughout we assume $0<q<1$. ($q$)-shifted factorials are given by
\begin{align*}
&(a)_k:=a(a+1)\ldots(a+k-1),\qquad\qquad\quad (a)_0:=1,\qquad
(a_1,\ldots,a_r)_k:=(a_1)_k\ldots(a_r)_k\,;\\
&(a;q)_k:=(1-a)(1-aq)\ldots(1-aq^{k-1}),\quad (a;q)_0:=1,\quad
(a_1,\ldots,a_r;q)_k:=(a_1;q)_k\ldots(a_r;q)_k\,.
\end{align*}
For $n$ a nonnegative integer we have terminating ($q$-)hypergeometric series
\begin{align*}
\hyp rs{-n,a_2,\ldots,a_r}{b_1,\ldots,b_s}z&:=
\sum_{k=0}^n \frac{(-n)_k}{k!}\,
\frac{(a_2,\ldots,a_r)_k}{(b_1,\ldots,b_s)_k}\,z^k,\\
\qhyp rs{-n,a_2,\ldots,a_r}{b_1,\ldots,b_s}{q,z}&:=
\sum_{k=0}^n \frac{(q^{-n};q)_k}{(q;q)_k}\,\frac{(a_2,\ldots,a_r;q)_k}
{(b_1,\ldots,b_s);q_k}\,\big((-1)^k q^{\half k(k-1)}\big)^{r-s+1}\,z^k.
\end{align*}
\section{Prelude: the one-variable case}
Let us explicitly consider the most simple situation, for polynomials
in one variable (in this section $n$ will denote the degree rather
than the number of variables).
Then both Jack and Macdonald polynomials are simple monomials $x^n$. Put
\begin{equation*}
P_n(x):=x^n,\quad P_n(x;q):=x^n,\quad
P_k^{\rm ip}(x):=x(x-1)\ldots(x-k+1)=(-1)^k(-x)_k.
\end{equation*}
$P_k^{\rm ip}(x)$ is the unique monic polynomial of degree $k$ which
vanishes at $0,1,\ldots,k-1$.
A binomial formula is given by
\begin{equation}
(x+1)^n=\sum_{k=0}^n \binom nk x^k,\quad{\rm or}\quad
P_n(x+1)=\sum_{k=0}^n \frac{P_k^{\rm ip}(n)}{P_k^{\rm ip}(k)}\,P_k(x).
\label{78}
\end{equation}
In the $q$-case put
\begin{equation*}
P_k^{\rm ip}(x;q):=(x-1)(x-q)\ldots(x-q^{k-1})=x^k(x^{-1};q)_k.
\end{equation*}
$P_k^{\rm ip}(x;q)$ is the unique monic polynomial of degree $k$ which
vanishes at
$1,q,\ldots,q^{k-1}$.
A $q$-binomial formula (see \cite[Exercise 1.6(iii)]{29}) is given by
\begin{align}
x^n=\qhyp20{q^{-n},x^{-1}}-{q,q^nx}&=
\sum_{k=0}^n\frac{(q^{-n},x^{-1};q)_k}{(-1)^k
q^{\half k(k-1)}(q;q)_k}\,(q^nx)^k,\quad{\rm or}
\nonumber\\
P_n(x;q)&=
\sum_{k=0}^n \frac{P_k^{\rm ip}(q^n;q)}{P_k^{\rm ip}(q^k;q)}\,
P_k^{\rm ip}(x;q).
\label{73}
\end{align}
Identity \eqref{78} is the limit case for $q\uparrow1$ of \eqref{73}.
The polynomials $P_k^{\rm ip}(x)$ and $P_k^{\rm ip}(x;q)$ are the
one-variable cases of the interpolation Jack and the interpolation
Macdonald polynomials, respectively.

In the one-variable case $BC_n$-type Jacobi polynomials become
classical Jacobi polynomials and Koornwinder polynomials become
Askey-Wilson polynomials. Their standard expressions as
($q$-)hypergeometric series are:
\begin{equation}
\frac{P_n^{(\al,\be)}(1-2x)}{P_n^{(\al,\be)}(1)}=
\sum_{k=0}^n \frac{(-n)_k(n+\al+\be+1)_k}{(\al+1)_k\,k!}\,x^k
=\hyp21{-n,n+\al+\be+1}{\al+1}x
\label{74}
\end{equation}
and
\begin{align}
\frac{p_n(\thalf(x+x^{-1});a_1,a_2,a_3,a_4\mid q)}
{p_n(\thalf(a_1+a_1^{-1});a_1,a_2,a_3,a_4\mid q)}
&=\sum_{k=0}^n
\frac{(q^{-n},q^{n-1}a_1a_2a_3a_4,a_1x,a_1x^{-1};q)_k}
{(a_1a_2,a_1a_3,a_1a_4,q;q)_k}\,q^k
\nonumber\\
&=\qhyp43{q^{-n},q^{n-1}a_1a_2a_3a_4,a_1x,a_1x^{-1}}
{a_1a_2,a_1a_3,a_1a_4}{q,q}.
\label{75}
\end{align}
Note that \eqref{74} gives an expansion in terms of monomials
$P_k(x)=x^k$ (Jack polynomials in one variable), while \eqref{75}
gives an expansion in terms of monic symmetric Laurent polynomials
\begin{equation}
P_k^{\rm ip}(x;q,a_1):=
\prod_{j=0}^{k-1}(x+x^{-1}-a_1q^j-a_1^{-1}q^{-j})
=\frac{(a_1x,a_1x^{-1};q)_k}{(-1)^k\,q^{\half k(k-1)}\,a_1^k}\,.
\label{76}
\end{equation}
The monic symmetric Laurent polynomial \eqref{76} is characterized by
its vanishing at $a_1$, $a_1q$, $\ldots$, $a_1q^{k-1}$. It is the
one-variable case of Okounkov's $BC$-type interpolation Macdonald
polynomial.  If we consider \eqref{75} as an expansion of its \LHS\ as
a function of $n$ then we see that it is expanded in terms of
functions $P_k^{\rm ip}(q^na_1';q,a_1')$ (using the definition in
\eqref{76}), where
$a_1':=(q^{-1}a_1a_2a_3a_4)^\half$.
Furthermore, if we replace $x$ by $a_1x$ in \eqref{76}, divide by
$a_1^k$, and let $a_1\to\iy$ then we obtain the $q$-binomial formula
\eqref{73}. Therefore, Okounkov \cite{1} calls \eqref{76}, as well as
its multi-variable analogue, also a binomial formula.

If we replace in \eqref{75} $a_1,a_2,a_3,a_4$ by
$q^{\al+1},-q^{\be+1},1,-1$ and let $q\uparrow1$ then we arrive at
\eqref{74}, which therefore might also be called a binomial formula.
If we consider \eqref{74} as an expansion of its \LHS\ as a
function of $n$ then we see
that it is expanded in terms of functions $P_k^{\rm ip}(n+\al';\al')$,
where $\al':=\thalf(\al+\be+1)$ and
\begin{equation}
P_k^{\rm ip}(x;\al):=\prod_{j=0}^{k-1}\big(x^2-(\al+j)\big)^2
=(-1)^k\,(\al-x)_k\,(\al+x)_k\,,
\end{equation}
a monic even polynomial of degree $2k$ in $x$ which is characterized
by its vanishing at $\al$, $\al+1$, $\ldots$, $\al+k-1$.
This is the one-variable case of the $BC$-type interpolation Jack
polynomial, which (possibly for the first time) will be defined in the
present paper.
\section{Preliminaries}
\subsection{Partitions}
\label{14}
We recapitulate some notions about partitions, diagrams and tableaux from
Macdonald \cite[\S I.1]{5}. However, different from \cite{5}, we
fix an integer $n\ge1$ and always understand
a {\em partition} $\la$ to be
of length $\le n$, i.e., $\la=(\la_1,\ldots,\la_n)\in\ZZ^n$ with
$\la_1\ge\la_2\ge\ldots\ge\la_n\ge0$.
Write $\ell(\la):=|\{j\mid \la_j>0\}|$ for the {\em length} of $\la$ and
$|\la|:=\la_1+\cdots+\la_n$ for its {\em weight}.
Also put
\begin{equation}
n(\la):=\sum_{i=1}^n(i-1)\la_i.
\label{44}
\end{equation}
We may abbreviate $k$ parts of $\la$ equal to $m$ by $m^k$ and we may
omit $0^k$ at the end.
For instance, $(2,2,1,0,0,0)=(2^2,1,0^3)=(2^2,1)$.
There is the special partition
\begin{equation}
\de:=(n-1,n-2,\ldots,1,0).
\label{45}
\end{equation}

A partition $\la$ can be displayed by a {\em Young diagram}, also
notated by $\la$, which consists of boxes
$(i,j)$ with $i=1,\ldots,\ell(\la)$ and $j=1,\ldots,\la_i$ for
a given $i$.
The {\em conjugate} partition $\la'$ has diagram such that
$(j,i)\in\la'$ iff $(i,j)\in\la$
The example below of the diagram of $\la=(7,5,5,2,2)$
and its conjugate $\la'=(5,5,3,3,3,1,1)$
will make clear how a diagram is drawn:
{\footnotesize\begin{center}
\yng(7,5,5,2,2)\qquad\qquad
\yng(5,5,3,3,3,1,1)
\end{center}}
For $(i,j)$ a box of a partition $\la$, the {\em arm-length} $a_\la(i,j)$
and {\em leg-length} $l_\la(i,j)$ are defined by
\begin{equation*}
a_\la(i,j):=\la_i-j,\qquad
l_\la(i,j):=|\{k>i\mid \la_k\ge j\}|.
\end{equation*}
Also the {\em arm-colength} $a_\la'(i,j)$
and {\em leg-colength} $l_\la'(i,j)$ are defined by
\begin{equation*}
a_\la'(i,j):=j-1,\qquad
l_\la'(i,j):=i-1.
\end{equation*}

The {\em dominance partial ordering} $\le$ and the
{\em inclusion partial ordering} $\subseteq$ are defined by:
\begin{align*}
\mu\le\la\quad&{\rm iff}\quad\mu_1+\cdots+\mu_i\le\la_1+\cdots+\la_i\quad
(i=1,\ldots,n);\\
\mu\subseteq\la\quad&{\rm iff}\quad\mu_i\le\la_i\quad(i=1,\ldots,n).
\end{align*}
Clearly, if $\mu\subseteq\la$ then $\mu\le\la$, while $\mu<\la$ implies that $\mu$ is less than $\la$
in the lexicographic ordering.
If $\mu\subseteq\la$
then we say that $\la$ {\em contains} $\mu$.
Note that, for the dominance partial ordering, we do not make the
usual requirement that $|\la|=|\mu|$.

For $\mu\subseteq\la$ define the {\em skew diagram} $\la-\mu$ as the
set of boxes
$\{s\in\la\mid s\notin\mu\}$.
A {\em horizontal strip} is a skew diagram with at most one box in each column.

For a horizontal strip $\la-\mu$ define $(R\backslash C)_{\la/\mu}$
as the set of boxes which are in a row of $\la$ intersecting with $\la-\mu$
but not in a column of $\la$ intersecting with $\la-\mu$. Then clearly
$(R\backslash C)_{\la/\mu}$ is completely contained in $\mu$.
For an example consider again $\la=(7,5,5,2,2)$ and take
$\mu=(5,5,3,2,1)$. In the following diagram
the cells of $\la-\mu$ have black squares and the cells of
$(R\backslash C)_{\la/\mu}$ have black diamonds.
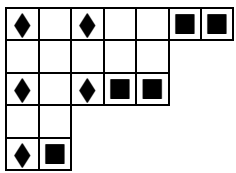
\begin{center}
\young(\blacklozenge\hfil\blacklozenge\hfil\hfil\blacksquare\blacksquare,\hfil\hfil\hfil\hfil\hfil,\blacklozenge\hfil\blacklozenge\blacksquare\blacksquare,\hfil\hfil,\blacklozenge\blacksquare)
\end{center}

\subsection{Tableaux}
\label{30}
For $\la$ a partition (of length $\le n$) we can fill the boxes $s$ of
$\la$ by numbers $T(s)\in\{1,2,\ldots,n\}$. Then $T$ is called a
{\em reverse tableau of shape $\la$ with entries in $\{1,\ldots,n\}$} if
$T(i,j)$ is weakly decreasing in $j$ and strongly decreasing in $i$.
(Clearly, the number of different entries has to be $\ge\ell(\la)$. In
\cite[\S I.1]{5} tableaux rather than reverse tableaux are defined.)\;
For an example consider again $\la=(7,5,5,2,2)$, so $\ell(\la)=5$. Let
us take $n=6$.
Then an example of a reverse tableau $T$ of shape $\la$ is given by
\begin{center}
\young(6664311,55522,44211,32,21)
\end{center}

For $T$ of shape $\la$ and for $k=0,1,\ldots,n$ let $\la^{(k)}$ be the
partition of which the
Young tableau consists of all $s\in\la$ such that
$T(s)>k$. Thus
\begin{equation}
0^n=\la^{(n)}\subseteq\la^{n-1}\subseteq\ldots\subseteq\la^{(1)}
\subseteq\la^{(0)}
=\la.
\label{5}
\end{equation}
Then the skew diagram $\la^{(k-1)}-\la^{(k)}$ is actually a horizontal strip
and it consists
of all boxes $s$ with $T(s)=k$.
We call the sequence $(\mu_1,\ldots,\mu_n)$ with
$\mu_k:=|\la^{(k-1)}-\la^{(k}|=|T^{-1}(\{k\})$ the {\em weight} of $T$.
For $\la$ and $T$ as in the example the inclusion sequence \eqref{5} becomes:
\[
()\subseteq(3)\subseteq(3,3)\subseteq(4,3,2)
\subseteq(5,3,2,1)\subseteq(5,5,3,2,1)\subseteq(7,5,5,2,2).
\]
and $T$ has weight $(5,5,2,3,3,3)$.

For a skew diagram $\la-\mu$ a {\em standard tableau} $T$ of shape $\la-\mu$
puts $T(s)$ in box $s$ of $\la-\mu$ such that each number in
$\{1,\ldots,|\la-\mu|\}$ occurs and $T(s)$ is strictly increasing in each row
and in each column.
\subsection{Symmetrized monomials}
Write $x^\mu:=x_1^{\mu_1}\ldots x_n^{\mu_n}$ for $\mu\in\ZZ^n$.
We say that $x^\mu$ has degree $|\mu|:=\mu_1+\cdots+\mu_n$. By the
degree of a Laurent polynomial $p(x)$ we mean the highest degree of a
monomial occurring in the Laurent expansion of $p(x)$.

Let $S_n$ be the symmetric group in $n$ letters and $W_n:=S_n\ltimes(\ZZ_2)^n$.
For $\la$ a partition and $x=(x_1,\ldots,x_n)\in\CC^n$ put
\begin{equation}
m_\la(x):=\sum_{\mu\in S_n\la} x^\mu,\qquad
\wt m_\la(x):=\sum_{\mu\in W_n\la} x^\mu.
\label{16}
\end{equation}
They form a basis of the space of $S_n$-invariant polynomials (respectively,
$W_n$-invariant Laurent polynomials) in $x_1,\ldots,x_n$.
Call an $S_n$-invariant polynomial (respectively,
$W_n$-invariant Laurent polynomial) of degree $|\la|$
{\em $\la$-monic} if its coefficient of
$m_\la$ (respectively, of $\wt m_\la$) is equal to 1.

\section{Macdonald and Koornwinder polynomials and $q=1$ limits}
From now on $n$ will be the number of variables and we will assume $n>1$.
\subsection{Macdonald polynomials}
See (4.7), (9.3), (9.5) and Remark on p.372 in Ch.~VI in
Macdonald \cite{5}.

Let $0<t<1$.
If $x=(x_1,\ldots,x_n)$ with $x_j\ne0$ for all $j$
then write $x^{-1}:=(x_1^{-1},\ldots,x_n^{-1})$.
Put
\[
\De_+(x)=\De_+(x;q,t):=\prod_{1\le i<j\le n}\frac{(x_i^{\,}x_j^{-1};q)_\iy}
{(tx_i^{\,}x_j^{-1};q)_\iy}\,,\qquad\De(x):=\De_+(x)\De_+(x^{-1}).
\]
{\em Macdonald polynomials} (for root system $A_{n-1}$) are $\la$-monic
$S_n$-invariant polynomials
\begin{equation}
P_\la(x;q,t)=P_\la(x)=
\sum_{\mu\le\la}u_{\la,\mu} m_\mu(x)
\label{10}
\end{equation}
such that (with $q,t$-dependence of $P_\la$ and $\De$ understood)
\begin{equation}
\int_{\TT^n} P_\la(x)\,m_\mu(x^{-1})\,
\De(x)\,\frac{dx_1}{x_1}\ldots\frac{dx_n}{x_n}=0\quad\mbox{if $\mu<\la$.}
\label{2}
\end{equation}
Here $\TT^n$ is the $n$-torus in $\CC^n$.
It follows from \eqref{2} that
\begin{equation}
\int_{\TT^n} P_\la(x)\,P_\mu(x^{-1})\,
\De(x)\,\frac{dx_1}{x_1}\ldots\frac{dx_n}{x_n}=0
\label{1}
\end{equation}
if $\mu<\la$, and that $P_\la$ is homogeneous of degree $|\la|$.
In fact, it can be shown that the orthogonality \eqref{1}
holds for $\la\ne\mu$. This deeper and very important result will also be met
for the three other orthogonal families discussed below.

Macdonald polynomials can be explicitly evaluated in a special point
(see \cite[Ch.~VI, (6.11)]{5}):
\begin{equation}
P_\la(t^\de;q,t)=t^{n(\la)}\,\frac{\De_+(q^\la t^\de;q,t)}{\De_+(t^\de;q,t)}
=t^{n(\la)}\,
\prod_{1\le i<j\le n}\frac{(t^{j-i+1};q)_{\la_i-\la_j}}
{(t^{j-i};q)_{\la_i-\la_j}}\,.
\label{46}
\end{equation}
There is also the duality result (see \cite[Ch.~VI, (6.6)]{5}):
\begin{equation}
\frac{P_\la(q^\nu t^\de;q,t)}{P_\la( t^\de;q,t)}
=\frac{P_\nu(q^\la t^\de;q,t)}{P_\nu(t^\de;q,t)}\,.
\label{55}
\end{equation}
\subsection{Jack polynomials}
See (10.13), (10.14), (10.35) and (10.36) in Ch.~VI in
Macdonald \cite{5} and see Stanley \cite{16}.

Let $\tau>0$. Put
\begin{equation}
\De_+(x)=\De_+(x;\tau):=\prod_{1\le i<j\le n}(1-x_i^{\,}x_j^{-1})^\tau,\qquad
\De(x):=\De_+(x)\De_+(x^{-1}).
\label{80}
\end{equation}
{\em Jack polynomials} are $\la$-monic
$S_n$-invariant polynomials
\[
P_\la(x;\tau)=P_\la(x)=
\sum_{\mu\le\la}u_{\la,\mu} m_\mu(x)
\]
satisfying \eqref{2} with $\De$ given by \eqref{80}.
Hence they satisfy \eqref{1}
if $\mu<\la$, and $P_\la$ is homogeneous of degree $|\la|$.
In fact, it can be shown that they satisfy \eqref{1} for $\la\ne\mu$.

Jack polynomials are limits of Macdonald polynomials:
\begin{equation}
\lim_{q\uparrow1}P_\la(x;q,q^\tau)=P_\la(x;\tau).
\label{22}
\end{equation}
Our notation of Jack polynomials relates to Macdonald's notation by
$P_\la(x;\tau)=P_\la^{(\tau^{-1})}(x)$.
Alternatively, \cite[Ch.~VI, (10.22)]{5} and \cite[Theorem 1.1]{16} work with
$J_\la^{(\al)}(x)=J_\la(x;\al)$, respectively. Then ((10.22) and (10.21)
in \cite[Ch.~VI]{5}):
\begin{equation*}
J_\la^{(\al)}=c_\la(\al)\,P_\la^{(\al)},\quad
c_\la(\al)=\prod_{s\in\la}(\al a(s)+l(s)+1).
\end{equation*}
We have the evaluation (see \cite[Theorem 5.4]{16}):
\begin{equation}
P_\la(1^n;\tau)=
\prod_{1\le i<j\le n}
\frac{((j-i+1)\tau)_{\la_i-\la_j}}{((j-i)\tau)_{\la_i-\la_j}}\,.
\label{47}
\end{equation}
The following limit is formally suggested by \eqref{22}:
\begin{equation}
\lim_{q\uparrow1} P_\la(q^{\tau\de};q,q^\tau)=P_\la(1^n;\tau).
\label{65}
\end{equation}
It follows rigorously by comparing \eqref{46} and \eqref{47}.
\subsection{Koornwinder polynomials}
See Koornwinder \cite{6}.

Let $|a_1|,|a_2|,|a_3|,|a_4|\le1$ such that $a_ia_j\ne1$ if $i\ne j$, and
such that non-real $a_j$
occur in complex conjugate pairs. Let $0<t<1$.
Put
\[
\De_+(x)=\De_+(x;q,t;a_1,a_2,a_3,a_4):=
\prod_{j=1}^n\frac{(x_j^2;q)_\iy}
{(a_1x_j,a_2x_j,a_3x_j,a_4x_j;q)_\iy}\,
\prod_{1\le i<j\le n}\frac{(x_ix_j,x_ix_j^{-1};q)_\iy}
{(tx_ix_j,tx_ix_j^{-1};q)_\iy}.
\]
Put
$\De(x):=\De_+(x)\De_+(x^{-1})$.
{\em Koornwinder polynomials} are $\la$-monic
$W_n$-invariant Laurent polynomials
\begin{equation}
P_\la(x;q,t;a_1,a_2,a_3,a_4)=P_\la(x)=
\sum_{\mu\le\la}u_{\la,\mu} \wt m_\mu(x)
\label{52}
\end{equation}
such that
\begin{equation}
\int_{\TT^n} P_\la(x)\,\wt m_\mu(x)\,
\De(x)\,\frac{dx_1}{x_1}\ldots\frac{dx_n}{x_n}=0\quad\mbox{if $\mu<\la$.}
\label{3}
\end{equation}
It follows from \eqref{3} that
\begin{equation}
\int_{\TT^n} P_\la(x)\,P_\mu(x)\,
\De(x)\,\frac{dx_1}{x_1}\ldots\frac{dx_n}{x_n}=0
\label{4}
\end{equation}
if $\mu<\la$, and that $P_\la$ is symmetric in $a_1,a_2,a_3,a_4$.
In fact, it can be shown that
the orthogonality \eqref{4} holds for $\la\ne\mu$.
Koornwinder polynomials
are a 5-parameter generalization of Macdonald's \cite{21}
3-parameter $q$-polynomials associated with root system $BC_n$.

Van Diejen \cite[\S5.2]{7} showed
that the Macdonald polynomial $P_\la(x;q,t)$ is the term
of highest degree $|\la|$ of $P_\la(x;q,t;a_1,a_2,a_3,a_4)$:
\begin{equation}
\lim_{r\to\iy} r^{-|\la|}P_\la(rx;q,t;a_1,a_2,a_3,a_4)=P_\la(x;q,t).
\label{8}
\end{equation}

For the following result we will need dual parameters $a_1',a_2',a_3',a_4'$:
\begin{equation}
a_1':=(q^{-1}a_1a_2a_3a_4)^\half,\quad
a_1'a_2'=a_1a_2,\quad a_1'a_3'=a_1a_3,\quad a_1'a_4'=a_1a_4.
\label{42}
\end{equation}
Below the ambiguity in taking a square root will cause no harm
because $\De_+$ is invariant under the transformation
$(x,a_1,a_2,a_3,a_4)\to(-x_1,-a_1,-a_2,-a_3,-a_4)$, by which $P_\la$
will also have this invariance, up to a factor $(-1)^{|\la|}$.

An evaluation formula for Koornwinder polynomials was conjectured by Macdonald
(1991, unpublished; see \cite[(5.5)]{17}). It reads:
\begin{align}
P_\la(t^\de a_1;q,t;a_1,a_2,a_3,a_4)
&=t^{-\lan\la,\de\ran}a_1^{-|\la|}\,
\frac{\De_+(q^\la t^\de a_1';q,t;a_1',a_2',a_3',a_4')}
{\De_+(t^\de a_1';q,t;a_1',a_2',a_3',a_4')}
\nonumber\\
&=t^{-\lan\la,\de\ran}a_1^{-|\la|}
\prod_{j=1}^n 
\frac{(t^{n-j}a_1'^2;q)_{\la_j}}{(t^{2n-2j}a_1'^2;q)_{2\la_j}}\,
(t^{n-j}a_1a_2,t^{n-j}a_1a_3,t^{n-j}a_1a_4;q)_{\la_j}
\nonumber\\
&\qquad\qquad\prod_{1\le i<j\le n}
\frac{(t^{2n-i-j+1} a_1'^2;q)_{\la_i+\la_j}}{(t^{2n-i-j}
a_1'^2;q)_{\la_i+\la_j}}\,
\frac{(t^{j-i+1};q)_{\la_i-\la_j}}{(t^{j-i};q)_{\la_i-\la_j}}\,.
\label{48}
\end{align}
It was proved by van Diejen \cite[(5.5)]{17} in the self-dual case $a_1=a_1'$.
In that case he also proved \cite[(5.4)]{17} Macdonald's duality
conjecture (1991):
\begin{equation}
\frac{P_\la(q^\nu t^\de a_1;q,t;a_1,a_2,a_3,a_4)}
{P_\la(t^\de a_1;q,t;a_1,a_2,a_3,a_4)}
=\frac{P_\nu(q^\la t^\de a_1';q,t;a_1',a_2',a_3',a_4')}
{P_\nu(t^\de a_1';q,t;a_1',a_2',a_3',a_4')}\,.
\label{49}
\end{equation}
Sahi \cite{18} proved \eqref{49} in the general case.
As pointed out in \cite[Section 7.2]{17},
this also implies \eqref{48} in the general case.
Macdonald independently proved his conjectures in his book \cite{20}, see there
(5.3.12) and (5.3.5), respectively.
\subsection{$BC_n$-type Jacobi polynomials}
See \cite{9} and \cite[Definition 3.5 and (3.18)]{8}.

Let $\al,\be>-1$ and $\tau>0$. Put
\[
\De(x)=\De(x;\tau;\al,\be):=
\prod_{j=1}^n x_j^\al(1-x_j)^\be\prod_{1\le i<j\le n}|x_i-x_j|^{2\tau}.
\]
{\em $BC_n$-type Jacobi polynomials}
are $\la$-monic $S_n$-invariant polynomials
\[
P_\la(x;\tau;\al,\be)=P_\la(x)=\sum_{\mu\le\la}a_{\la,\mu} m_\mu(x)
\]
such that
\begin{equation}
\int_{[0,1]^n} P_\la(x)\,m_\mu(x)\,\De(x)\,dx_1\ldots dx_n=0\quad
\mbox{if $\mu<\la$.}
\label{7}
\end{equation}
It follows from \eqref{7} that
\begin{equation}
\int_{[0,1]^n} P_\la(x)\,P_\mu(x)\,\De(x)\,dx_1\ldots dx_n=0
\label{6}
\end{equation}
if $\mu<\la$. It can be shown, see \cite[Corollary 3.12]{9},
that \eqref{6} holds more generally if $\la\ne\mu$.

The case $c=1$, $d=-1$ of \cite[(5.5)]{8} says that
\begin{equation}
\lim_{q\uparrow1}P_\la(x;q,q^\tau;q^{\al+1},-q^{\be+1},1,-1)=
(-4)^{|\la|}\,P_\la(\tfrac14(2-x-x^{-1});\tau;\al,\be).
\label{58}
\end{equation}
Furthermore, it was pointed out in \cite[(4.8)]{10}
that the Jack polynomial $P_\la(x;\tau)$
is the term
of highest degree $|\la|$ of the
$BC_n$-type Jacobi polynomial $P_\la(x;\tau;\al,\be)$:
\begin{equation}
\lim_{r\to\iy} r^{-|\la|}P_\la(rx;\tau;\al,\be)=P_\la(x;\tau).
\label{9}
\end{equation}
This is the $q=1$ analogue of  the limit \eqref{8}.

An evaluation formula for Jacobi polynomials associated with root systems,
including $BC_n$,
was given by Opdam \cite[Corollary 5.2]{22}.
See reformulations of this result in the $BC_n$ case by
van Diejen \cite[(6.43d)]{23} and by Halln\"as \cite[p.1594]{24}.
The formula can be given very explicitly as follows:
\begin{multline}
 P_\la(0;\tau;\al,\be)=(-1)^{|\la|}\,\prod_{j=1}^n
\frac{((n-j)\tau+2\al')_{\la_j}\,
((n-j)\tau+\al+1)_{\la_j}}{((2n-2j)\tau+2\al')_{2\la_j}}\\
\times\prod_{1\le i<j\le n}
\frac{((2n-i-j+1)\tau+2\al')_{\la_i+\la_j}}
{((2n-i-j)\tau+2\al')_{\la_i+\la_j}}\,
\frac{((j-i+1)\tau)_{\la_i-\la_j}}{((j-i)\tau)_{\la_i-\la_j}}\,.
\label{59}
\end{multline}
Here
\begin{equation}
\al':=\thalf(\al+\be+1).
\label{81}
\end{equation}

The following limit is formally suggested by \eqref{58}:
\begin{equation}
\lim_{q\uparrow1}P_\la(q^{\tau\de+\al+1};q,q^\tau;q^{\al+1},-q^{\be+1},1,-1)=
(-4)^{|\la|}\,P_\la(0;\tau;\al,\be).
\label{60}
\end{equation}
It follows rigorously by comparing \eqref{48} and \eqref{59}.
\section{Interpolation polynomials}
\subsection{Interpolation Macdonald polynomials}
See Sahi \cite[Theorem 1.1]{12},  Knop \cite[Theorem 2.4(b)]{13} and
Okounkov \cite[(4.2), (4.3)]{11}.

Let $0<t<1$.
The {\em interpolation Macdonald polynomial}
(or {\em shifted Macdonald polynomial})
$P_\la^{\rm ip}(x;q,t)$ is the unique $\la$-monic $S_n$-invariant
polynomial of degree $|\la|$ such that
$P_\la^{\rm ip}(q^\mu t^\de;q,t)=0$
for each partition $\mu\ne\la$ having $|\mu|\le|\la|$.
Here $q^\mu t^\de=(q^{\mu_1}t^{n-1},q^{\mu_2}t^{n-2},\ldots,q^{\mu_n})$.

Our $P_\la^{\rm ip}$ is related to Sahi's $R_\la$, Knop's $P_\la$
(use \cite[Theorem 3.11]{13})
and
Okounkov's $P_\la^*$ (use \cite[(4.11)]{11}), respectively,  as follows:
\[
P_\la^{\rm ip}(x;q,t)=R_\la(x;q^{-1},t^{-1})=t^{(n-1)|\la|}\,P_\la(t^{-(n-1)}x)
=t^{(n-1)|\la|}\,P_\la^*(xt^{-\de}).
\]
Okounkov \cite{11} speaks about {\em shifted} polynomials because in
his notation the polynomials are only symmetric after a
(multiplicative) shift.

$P_\la^{\rm ip}$ has the extra vanishing property (see \cite[p.93]{13} or
\cite[(4.12)]{11}):
\[
P_\la^{\rm ip}(q^\mu t^\de;q,t)=0\quad\mbox{if $\mu$ is a partition
not containing $\la$}.
\]

By \cite[(4.11)]{11} $P_\la^{\rm ip}$ can be expanded in terms of
Macdonald polynomials as follows:
\begin{equation}
P_\la^{\rm ip}(x;q,t)=\sum_{\mu\subseteq\la}b_{\la,\mu}\,P_\mu(x;q,t)
\label{11}
\end{equation}
for certain coefficients $b_{\la,\mu}$, where $b_{\la,\la}=1$ by
$\la$-monicity.
This has several consequences.  By combination with \eqref{10} we see that
\begin{equation}
P_\la^{\rm ip}(x;q,t)=\sum_{\mu\le\la}c_{\la,\mu}\,m_\mu(x)
\label{12}
\end{equation}
for certain coefficients $c_{\la,\mu}$, and $c_{\la,\la}=1$.
Furthermore, by \eqref{11}, $P_\la(x;q,t)$ is the term of highest
degree $|\la|$ of the polynomial $P_\la^{\rm ip}(x;q,t)$:
\begin{equation}
\lim_{r\to\iy} r^{-|\la|}P_\la^{\rm ip}(rx;q,t))=P_\la(x;q,t).
\label{13}
\end{equation}
Although Knop \cite{13} did not give \eqref{11}, he did give
\eqref{12} and \eqref{13},
proved differently (see Theorem 3.11 and Theorem 3.9 in \cite{13}).
The result \eqref{13} is also proved by Sahi \cite[Theorem~1.1]{12}.
\subsection{Interpolation Jack polynomials}
See Sahi \cite[Theorem 1]{37},
Knop \& Sahi \cite[p.475, 478]{15},
Okounkov \& Olshanski \cite[p.70]{14} and Okounkov \cite[Section~7]{11}.

Let $\tau>0$.
The {\em interpolation Jack polynomial} (or {\em shifted Jack polynomial})
$P_\la^{\rm ip}(x;\tau)$ is the unique $\la$-monic $S_n$-invariant
polynomial of degree $|\la|$ such that
$P_\la^{\rm ip}(\mu+\tau\de;\tau)=0$
for each partition $\mu\ne\la$ having $|\mu|\le|\la|$.
It can be expressed in terms of $P_\la^{\tau\de}$ from \cite{15}
and in terms of
$P_\la^*(\,.\,;\tau)$ from \cite{14}, \cite{11} as follows:
\begin{equation*}
P_\la^{\rm ip}(x;\tau)=P_\la^{\tau\de}(x)=P_\la^*(x-\tau\de;\tau).
\end{equation*}
It has the extra vanishing property (see \cite[Theorem 5.2]{15}):
\begin{equation*}
P_\la^{\rm ip}(\mu+\tau\de;\tau)=0\quad\mbox{if $\mu$ is a partition
not containing $\la$}.
\end{equation*}
It has an expansion of the form
\begin{equation*}
P_\la^{\rm ip}(x;\tau)=\sum_{\mu\le\la}c_{\la,\mu}\,m_\mu(x)
\end{equation*}
for certain coefficients $c_{\la,\mu}$, and $c_{\la,\la}=1$
(see \cite[Corollary 4.6]{15}).
The term of highest degree $|\la|$ of the
polynomial $P_\la^{\rm ip}(x;\tau)$ is the
Jack polynomial $P_\la(x;\tau)$ see \cite[Corollary 4.7]{15}):
\begin{equation}
\lim_{r\to\iy} r^{-|\la|}P_\la^{\rm ip}(rx;\tau)=P_\la(x;\tau).
\label{27}
\end{equation}

Interpolation Jack polynomials are a limit case of interpolation
Macdonald polynomials
(see \cite[(7.1)]{11}):
\begin{equation}
\lim_{q\uparrow1}(q-1)^{-|\la|} P_\la^{\rm ip}(q^x;q,q^\tau)=
P_\la^{\rm ip}(x;\tau).
\label{25}
\end{equation}
\subsection{$BC_n$-type interpolation Macdonald polynomials}
\begin{definition}[{\em $BC_n$-type interpolation}
(or {\em $BC_n$-type shifted}) {\em Macdonald polynomials}]
\quad\\
Let $0<t<1$ and  let $a\in\CC$ be generic.
$P_\la^{\rm ip}(x;q,t,a)$ is the unique $W_n$-invariant $\la$-monic Laurent
polynomial of degree $|\la|$ such that
\begin{equation}
\label{35}
P_\la^{\rm ip}(q^\mu t^\de a;q,t,a)=0\quad\mbox{if $\mu$ does not
contain $\la$.}
\end{equation}
(In particular, $P_\la^{\rm ip}(q^\mu t^\de a;q,t,a)=0$
if $|\mu|\le|\la|$, $\mu\ne\la$.)
\end{definition}

These polynomials were first introduced by Okounkov \cite{1} in a
different notation and
normalization. Okounkov \cite[p.185, Section 1]{1} specifies the
genericity of the parameter $a$ (in his notation $s$) by the condition
$q^i t^j a^k\ne1$ for $i,j,k\in\ZZ_{>0}$. However, this may be too
strong on the one hand and too weak on the other hand.
Certainly the \RHS\ of \eqref{104} (equivalently the normalization
constant in \cite[Definition 1.2]{1}) should be nonzero. A requirement
for this is that $q^i t^j a^2\ne1$ for $i\in\ZZ_{>0}$,
$j\in\ZZ_{\ge0}$.

Different approaches were given by Rains \cite{2},
and later by Noumi in unpublished slides of a lecture given in 2013 at
a conference at Kyushu University.
Our normalization follows Rains \cite{2}.
In terms of Rains' $\bar P_\la^{*(n)}$ and Okounkov's $P_\la^*$ we have
(cf.~\cite[p.76, Remark 1]{2}):
\begin{equation}
P_\la^{\rm ip}(x;q,t,a)=\bar P_\la^{*(n)}(x;q,t,a)
=(t^{n-1}a)^{|\la|}P_\la^*(xt^{-\de} a^{-1};q,t,a).
\label{50}
\end{equation}
Note that Okounkov's $P_\la^*(x;q,t,s)$ is $W_n$-symmetric in the variables
$x_it^{n-i}s$ ($i=1,2,\ldots,n$).

Just as for $P_\la^{\rm ip}(x;q,t)$, the top homogeneous term of
$P_\la^{\rm ip}(x;q,t,a)$ equals the Macdonald polynomial
$P_\la(x;q,t)$
(see \cite[Section 4]{1}):
\begin{equation}
\lim_{r\to\iy} r^{-|\la|}P_\la^{\rm ip}(rx;q,t,a))=P_\la(x;q,t).
\label{20}
\end{equation}
There is also a limit from $P_\la^{\rm ip}(x;q,t,a)$ to $P_\la^{\rm ip}(x;q,t)$
(see \cite[p.75]{2}):
\begin{equation}
\lim_{a\to\iy} a^{-|\la|}P_\la^{\rm ip}(ax;q,t,a)=P_\la^{\rm ip}(x;q,t).
\label{21}
\end{equation}

By combination of \eqref{50} with \cite[Definition 1.1 and 1.2]{1} we get the
following evaluation formula (with $\lan\,.\,,\,.\,\ran$ denoting the
standard inner product on $\RR^n$):
\begin{equation}
P_\la^{\rm ip}(q^\la t^\de a;q,t,a)=
q^{-\lan\la,\la\ran} t^{-\lan\la,\de\ran}\,a^{-|\la|}
\prod_{(i,j)\in\la}(1-q^{\la_i-j+1}t^{\la_j'-i})
(1-a^2 q^{\la_i+j-1}t^{\la_j'-i+2(n-\la_j')}).
\label{104}
\end{equation}
By \cite[Lemma 2.1]{2} (see also \cite[Corollary 3.7]{2})
this can be rewritten as
\begin{align}
P_\la^{\rm ip}(q^\la t^\de a;q,t,a)&=
q^{-\lan\la,\la\ran} t^{-\lan\la,\de\ran}\,a^{-|\la|}
\prod_{j=1}^n \frac{(q t^{n-j};q)_{\la_j}
(t^{2n-2j} a^2;q)_{2\la_j}}{(t^{n-j} a^2;q)_{\la_j}}
\nonumber\\
&\qquad\times\prod_{1\le i<j\le n}
\frac{(t^{2n-i-j} a^2;q)_{\la_i+\la_j}}{(t^{2n-i-j+1} a^2;q)_{\la_i+\la_j}}\,
\frac{(qt^{j-i-1};q)_{\la_i-\la_j}}{(q t^{j-i};q)_{\la_i-\la_j}}\,.
\label{51}
\end{align}

An elementary consequence of the definition of $P_\la^{\rm ip}(x;q,t,a)$
is a reduction formula (see \cite[Proposition 2.1]{1}):
\begin{equation}
P_\mu^{\rm ip}(x;q,t,a)=(-a)^{n\mu_n} q^{-\half n\mu_n(\mu_n-1)}
\prod_{j=1}^n\big((x_ja;q)_{\mu_n}(x_j^{-1}a;q)_{\mu_n}\big)
P_{\mu-\mu_n 1^n}^*(x;q,t,q^{\mu_n}a).
\label{69}
\end{equation}
\section{Combinatorial formulas}
\subsection{Combinatorial formula for Macdonald polynomials}
The combinatorial formula for Macdonald polynomials which we will use
is a special case of
Macdonald \cite[Ch.~VI, (7.13${}'$)]{5}:
\begin{equation}
P_\la(x;q,t)=\sum_T \psi_T(q,t) \prod_{s\in\la} x_{T(s)},
\label{17}
\end{equation}
where the sum is over all tableaux $T$ of shape $\la$ with entries in
$\{1,\ldots,n\}$
and with $\psi_T(q,t)$ defined in \cite[Ch.~VI]{5} by formula (7.11${}'$),
by formula (ii) on p.341 with
$C_{\la/\mu}$ and $S_{\la/\mu}$ given after (6.22), and by formula (6.20).
See \cite[Section 1]{35} for a summary of these results.
Since the Macdonald polynomial is symmetric, we may as well sum over
reverse tableaux
instead of tableaux, with the definition of $\psi_T(q,t)$ accordingly adapted.
We will now give $\psi_T(q,t)$ explicitly.
See Subsections \ref{14} and \ref{30} for notation.

Recall that for a reverse tableaux $T$ of shape $\la$ with entries in
$\{1,\ldots,n\}$ we write $0^n=\la^{(n)}\subseteq
\la^{(n-1)}\subseteq\ldots\subseteq \la^{(0)}=\la$ with $T(s)=i$ for
$s$ in the horizontal strip $\la^{(i-1)}-\la^{(i)}$.
Now
\begin{equation}
\psi_T(q,t):=\prod_{i=1}^n \psi_{\la^{(i-1)}/\la^{(i)}}(q,t),\qquad
\psi_{\mu/\nu}(q,t)=
\prod_{s\in (R\backslash C)_{\mu/\nu}}
\frac{b_\nu(s;q,t)}{b_\mu(s;q,t)}\,,
\label{24}
\end{equation}
where
\begin{equation}
b_\mu(s;q,t):=\frac{1-q^{a_\mu(s)} t^{l_\mu(s)+1}}{1-q^{a_\mu(s)+1}
t^{l_\mu(s)}}\,.
\label{82}
\end{equation}

By \eqref{16} and \eqref{10} we obtain that
\begin{equation}
P_\la(x;q,t)=\sum_{\mu\le\la}u_{\la,\mu}(q,t) m_\mu(x)\quad
\mbox{with}\quad
u_{\la,\mu}(q,t)=\sum_T \psi_T(q,t),
\label{23}
\end{equation}
where the $T$-sum is over all reverse tableaux $T$ of shape $\la$ and
weight $\mu$,
see \cite[p.378]{5}.

Let $T_\la$ be
the tableau of shape $\la$ for which $T(i,j)=n+1-i$. This is the
unique tableau of shape $\la$ which has weight
$(\la_n,\la_{n-1},\ldots,\la_1)$. Since $P_\la(x;q,t)$
is $\la$-monic, it follows from \eqref{17} that
\begin{equation}
\psi_{T_\la}(q,t)=1.
\label{33}
\end{equation}
\subsection{Combinatorial formulas
for ($BC_n$) interpolation Macdonald polynomials}
For interpolation Macdonald polynomials $P_\la^{\rm ip}(x;q,t)$
and $BC_n$-type interpolation Macdonald polynomials $P_\la^{\rm ip}(x;q,t,a)$
there are
combinatorial formulas similar to \eqref{17} and also involving
$\psi_T(q,t)$ given by
\eqref{24}:
\begin{align}
P_\la^{\rm ip}(x;q,t)&=\sum_T \psi_T(q,t) \prod_{s\in\la}
\big(x_{T(s)}-q^{a_\la'(s)} t^{n-T(s)-l_\la'(s)}\big),
\label{18}
\\
P_\la^{\rm ip}(x;q,t,a)&=\sum_T \psi_T(q,t) \prod_{s\in\la}
\Big(x_{T(s)}-q^{a_\la'(s)} t^{n-T(s)-l_\la'(s)}a\Big)
\,\Big(1-\big(q^{a_\la'(s)} t^{n-T(s)-l_\la'(s)}a\big)^{-1}x_{T(s)}^{-1}\Big).
\label{19}
\end{align}
See \cite[(1.4)]{11} for \eqref{18} and \cite[(5.3)]{1} for \eqref{19}.
The sums are over all reverse tableaux $T$ of shape $\la$ with entries in
$\{1,\ldots,n\}$.
Note that the limits \eqref{13}, \eqref{20} and \eqref{21} also follow
by comparing \eqref{17}, \eqref{18} and \eqref{19}.
\subsection{Combinatorial formulas for (interpolation) Jack polynomials}
The combinatorial formula for Jack polynomials can be obtained as a limit
case of the combinatorial formula \eqref{17} for Macdonald polynomials
by using the limit \eqref{22} (see \cite[p.379]{5}),
but it can also be obtained independently,
as was first done by Stanley \cite[Theorem 6.3]{16}:
\begin{equation}
P_\la(x;\tau)=\sum_T \psi_T(\tau) \prod_{s\in\la} x_{T(s)},
\label{28}
\end{equation}
where the sum is over all reverse tableaux $T$ of shape $\la$ with entries in
$\{1,\ldots,n\}$. For the definition of $\psi_T(\tau)$ take
$0^n=\la^{(n)}\subseteq \la^{(n-1)}\subseteq\ldots\subseteq \la^{(0)}=\la$
as before and put
\begin{equation}
\psi_T(\tau):=\prod_{i=1}^n\,
\prod_{s\in (R\backslash C)_{\la^{(i-1)}/\la^{(i)}}}
\frac{b_{\la^{(i)}}(s;\tau)}{b_{\la^{(i-1)}}(s;\tau)}
\label{26}
\end{equation}
with
\begin{equation*}
b_\mu(s;\tau):=\frac{a_\mu(s)+\tau(l_\mu(s)+1)}{a_\mu(s)+\tau l_\mu(s)+1}\,.
\end{equation*}
Note that
\begin{equation*}
\lim_{q\uparrow1}b_\mu(s;q,q^\tau)=b_\mu(s,\tau)\quad\mbox{and}\quad
\lim_{q\uparrow1}\psi_T(q,q^\tau)=\psi_T(\tau).
\end{equation*}
Hence, by \eqref{33} we have
\begin{equation}
\psi_{T_\la}(\tau)=1.
\label{34}
\end{equation}

Similarly as for \eqref{23} we derive immediately that
\begin{equation*}
P_\la(x;\tau)=\sum_{\mu\le\la}u_{\la,\mu}(\tau) m_\mu(x)\quad
\mbox{with}\quad
u_{\la,\mu}(\tau)=\sum_T \psi_T(\tau),
\end{equation*}
where the $T$-sum is over all reverse tableaux $T$ of shape $\la$ and
weight $\mu$.

The combinatorial formula for interpolation Jack polynomials
(see \cite[(2.4)]{14}) was obtained in \cite[Section 7]{11} as a limit
case of the combinatorial formula \eqref{18} for interpolation Macdonald
polynomials by using the limit \eqref{25}:
\begin{equation}
P_\la^{\rm ip}(x;\tau)=\sum_T \psi_T(\tau) \prod_{s\in\la}
\big(x_{T(s)}-a_\la'(s)-\tau(n-T(s)-l_\la'(s))\big),
\label{29}
\end{equation}
with the sum over all tableaux $T$ of shape $\la$ with entries in
$\{1,\ldots,n\}$ and $\psi_T(\tau)$ given by \eqref{26}.
The limit \eqref{27} also follows by comparing \eqref{29} and \eqref{28}.
Furthermore, by comparing \eqref{19} and \eqref{28} we obtain the limit
\begin{equation}
\lim_{q\uparrow1}P_\la^{\rm ip}(x;q,q^\tau,q^\al)=
P_\la(x+x^{-1}-2;\tau),
\label{36}
\end{equation}
and from \eqref{18} and \eqref{28} we obtain
\begin{equation}
\lim_{q\uparrow1} P_\la^{\rm ip}(x;q,q^\tau)=P_\la(x-1^n;\tau).
\label{39}
\end{equation}
\begin{remark}
\label{103}
When we compare combinatorial formulas in the case of $n$ and of $n-1$
variables, we see that, in general, a combinatorial formula is equivalent
to a {\em branching formula}, which expands a polynomial $P_\la$ in
$x_1,\ldots,x_n$ in terms of polynomials $P_\mu$ in $x_1,\ldots,x_{n-1}$
with the expansion coefficients depending on $x_n$.
In particular, the branching formula for Macdonald polynomials is
(see \cite[(1.9), (1.8)]{35})
\begin{equation}
P_\la(x_1,\ldots,x_{n-1},x_n;q,t))=
\sum_\mu P_{\la/\mu}(x_n;q,t)\,P_\mu(x_1,\ldots,x_{n-1};q,t),
\label{101}
\end{equation}
where the sum runs over all partitions $\mu\subseteq\la$ of length $<n$
such that
$\la-\mu$ is a horizontal strip, and where, in notation \eqref{24},
\begin{equation}
P_{\la/\mu}(z;q,t)=\psi_{\la/\mu}(q,t)\,z^{|\la|-|\mu|}.
\label{102}
\end{equation}
The coefficients $\psi_{\la/\mu}(q,t)$ can be expressed in terms of
Pieri coefficients for Macdonald polynomials by interchanging $q$ and
$t$ and by passing to conjugate partitions $\la',\mu'$:
$\psi_{\la/\mu}(q,t)=\psi_{\la'/\mu'}'(t,q)$, see \cite[Ch.~VI, (6.24)]{5}.

Van Diejen \& Emsiz \cite{36} recently obtained a branching formula
for Koornwinder polynomials. It has the same structure as \eqref{101},
but the analogue of \eqref{102} becomes a sum of terms in the \RHS.
Each term is similar to the \RHS\ of \eqref{102}, with the monomial being
replaced by a quadratic $q$-factorial, also depending on the parameter $a_1$.
The analogues of the coefficients $\psi_{\la/\mu}$ can be expressed in
terms of (earlier known) Pieri-type coefficients for Koornwinder
polynomials. By taking highest degree parts in both sides of
the new branching formula and by using \eqref{8}, we are reduced to
\eqref{102}. However, the combinatorial formula \eqref{19} and its
corresponding branching formula for $BC_n$-type interpolation Macdonald
polynomials are quite different from the results in \cite{36}.
 
\end{remark}
\section{$BC_n$-type interpolation Jack polynomials}
In view of the results surveyed until now the following definition is
quite natural:
\begin{definition}
Let $\tau>0$ and let $\al\in\CC$ be generic.
The {\em $BC_n$-type interpolation Jack polynomial}
$P_\la^{\rm ip}(x;\tau,\al)$ is given as a limit
of $BC_n$-type interpolation Macdonald polynomials by:
\begin{equation}
P_\la^{\rm ip}(x;\tau,\al):=\lim_{q\uparrow1}
(1-q)^{-2|\la|} P_\la^{\rm ip}(q^x;q,q^\tau,q^\al).
\label{31}
\end{equation}
\end{definition}
Concerning the genericity of $\al\in\CC$ we should have at least
that the evaluation \eqref{63} is nonzero, i.e.,
$i+j\tau+2\al\ne0$ for $i\in\ZZ_{>0}$ and $j\in\ZZ_{\ge0}$.
That the limit \eqref{31}
exists can be seen by substituting \eqref{19} in the
\RHS\ of \eqref{31}. We obtain
\begin{equation}
P_\la^{\rm ip}(x;\tau,\al)=\sum_T \psi_T(\tau)\,
\prod_{s\in\la}
\Big(x_{T(s)}^2-\big(a_\la'(s)+\tau(n-T(s)-l_\la'(s))+\al\big)^2\Big)
\label{32}
\end{equation}
with the sum over all reverse tableaux $T$ of shape $\la$ with entries in
$\{1,\ldots,n\}$ and $\psi_T(\tau)$ given by \eqref{26}.
From \eqref{31}, \eqref{32} and the properties of
$P_\la^{\rm ip}(x;q,t)$
we see that $P_\la^{\rm ip}(x;\tau,\al)$ is a $W_n$-invariant polynomial
of degree $2|\la|$ in $x$, where $(\ZZ_2)^n$ now acts on the polynomial
by sending some of the variables $x_i$ to $-x_i$ rather than to $x_i^{-1}$.
By \eqref{34} it follows from \eqref{32} that $P_\la^{\rm ip}(x;\tau,\al)$
is ($2\la$)-monic.

It follows from \eqref{35} and \eqref{31} that
\begin{equation*}
P_\la^{\rm ip}(\mu+\tau\de+\al;\tau,\al)=0\quad\mbox{if $\mu$
does not contain $\la$.}
\end{equation*}
By comparing \eqref{32}, \eqref{29} and \eqref{28} we obtain the limits
\begin{align}
\lim_{r\to\iy}r^{-2|\la|}P_\la^{\rm ip}(rx;\tau,\al)&=P_\la(x^2;\tau),
\label{40}\\
\lim_{\al\to\iy}(2\al)^{-|\la|}P_\la^{\rm ip}(x+\al;\tau,\al)&=
P_\la^{\rm ip}(x;\tau).
\label{41}
\end{align}

From \eqref{104} and
\eqref{51} together with \eqref{31} we obtain the evaluation formula
\begin{align}
P_\la^{\rm ip}(\la+\tau\de+\al;\tau,\al)&=
\prod_{(i,j)\in\la} \big(\la_i-j+1+\tau(\la_j'-i)\big)
\big(2\al+\la_i+j-1+\tau(\la_j'-i+2(n-\la_j')\big)
\nonumber\\
&=\prod_{j=1}^n \frac{((n-j)\tau+1)_{\la_j} (2(n-j)\tau+2\al)_{2\la_j}}
{((n-j)\tau+2\al)_{\la_j}}
\nonumber\\
&\quad\times\prod_{1\le i<j\le n}
\frac{((2n-i-j)\tau+2\al)_{\la_i+\la_j}}
{((2n-i-j+1)\tau+2\al)_{\la_i+\la_j}}\,
\frac{((j-i-1)\tau+1)_{\la_i-\la_j}}{((j-i)\tau+1)_{\la_i-\la_j}}\,.
\label{63}
\end{align}
By \eqref{69} and \eqref{31} we get a reduction formula
\begin{equation}
P_\la^{\rm ip}(x;\tau,\al)=(-1)^{n\la_n}\,
\prod_{j=1}^n\big((\al+x_j)_{\la_n}(\al-x_j)_{\la_n}\big)\,
P_{\la-\la_n 1^n}^{\rm ip}(x;\tau,\la_n+\al).
\label{70}
\end{equation}
\section{Binomial formulas}
\subsection{Binomial formula for Koornwinder polynomials}
Okounkov \cite[Theorem 7.1]{1} obtained the {\em binomial formula} for
Koornwinder polynomials:\begin{equation}
\frac{P_\la(x;q,t;a_1,a_2,a_3,a_4)}{P_\la(t^\de a_1;q,t;a_1,a_2,a_3,a_4)}
=\sum_{\mu\subseteq\la}
\frac{P_\mu^{\rm ip}(q^\la t^\de a_1';q,t,a_1')}
{P_\mu^{\rm ip}(q^\mu t^\de a_1';q,t,a_1')}\,
\frac{P_\mu^{\rm ip}(x;q,t,a_1)}{P_\mu(t^\de a_1;q,t;a_1,a_2,a_3,a_4)}\,.
\label{43}
\end{equation}
As pointed out in \cite{1}, the duality \eqref{49} immediately follows
from \eqref{43} in the self-dual case $a_1=a_1'$. In the general case
\eqref{49} will follow from \eqref{43} together with the identity
\begin{equation*}
\frac{P_\mu(t^\de a_1;q,t;a_1,a_2,a_3,a_4)}
{P_\mu^{\rm ip}(q^\mu t^\de a_1;q,t,a_1)}=
\frac{P_\mu(t^\de a_1';q,t;a_1',a_2',a_3',a_4')}
{P_\mu^{\rm ip}(q^\mu t^\de a_1';q,t,a_1')}\,,
\end{equation*}
which is a consequence of the evaluation formulas \eqref{48} and \eqref{51}.

Rains \cite[Section 5]{2} gives an alternative definition of
Koornwinder polynomials involving triangularity and evaluation
symmetry by which a version of \eqref{43} is an immediate consequence.

Because the inclusion partial ordering is compatible with the
lexicographic ordering, we can use induction in $\la$ with respect to
the lexicographic ordering in order to show from \eqref{43} that
\begin{equation*}
P_\la^{\rm ip}(x;q,t,a_1)
=\sum_{\mu\subseteq\la} b_{\la,\mu}\,P_\mu(x;q,t;a_1,a_2,a_3,a_4)
\end{equation*}
for certain coefficients $b_{\la,\mu}$
(more explicitly given in \cite[Theorem 5.12]{2}).
Together with \eqref{52} this implies that
\begin{equation}
P_\la^{\rm ip}(x;q,t,a_1)=\sum_{\mu\le\la}c_{\la,\mu}\,\wt m_\mu(x)
\label{53}
\end{equation}
for certain  coefficients $c_{\la,\mu}$.
\subsection{Binomial formula for Macdonald polynomials}
Okounkov \cite[(1.11)]{19} gave the binomial formula for Macdonald polynomials:
\begin{equation}
\frac{P_\la(x;q,t)}{P_\la(t^\de;q,t)}
=\sum_{\mu\subseteq\la}
\frac{P_\mu^{\rm ip}(q^\la t^\de;q,t)}{P_\mu^{\rm ip}(q^\mu t^\de;q,t}\,
\frac{P_\mu^{\rm ip}(x;q,t)}{P_\mu(t^\de;q,t)}\,.
\label{54}
\end{equation}
The binomial formula \eqref{54} immediately implies
the duality formula \eqref{55}
for Macdonald polynomials.

By comparing the binomial formulas \eqref{43} and \eqref{54} we can
obtain a new limit formula:
\begin{theorem}
There is the limit
\begin{equation}
\lim_{a_1\to\iy} a_1^{-|\la|}\,P_\la(a_1x;q,t;a_1,a_2,a_3,a_4)=P_\la(x;q,t).
\label{56}
\end{equation}
\end{theorem}
\Proof
Rewrite \eqref{43} as
\[
\frac{a_1^{-|\la|}\, P_\la(a_1x;q,t;a_1,a_2,a_3,a_4)}
{a_1^{-|\la|}\,P_\la(t^\de a_1;q,t;a_1,a_2,a_3,a_4)}
=\sum_{\mu\subseteq\la}
\frac{a_1^{-|\mu|}P_\mu^{\rm ip}(q^\la t^\de a_1';q,t,a_1')}
{a_1^{-|\mu|}P_\mu^{\rm ip}(q^\mu t^\de a_1';q,t,a_1')}\,
\frac{a_1^{-|\mu|}\,P_\mu^{\rm ip}(a_1x;q,t,a_1)}{a_1^{-|\mu|}\,
P_\mu(t^\de a_1;q,t;a_1,a_2,a_3,a_4)}\,.
\]
Now let $a_1\to\iy$. The result follows by \eqref{21} and \eqref{54}
if we can show that \eqref{56} holds for $x=t^\de$. This, in its turn,
follows by comparing \eqref{48} and \eqref{46}.
\qed
\bPP
As pointed out to me by Ole Warnaar, the limit \eqref{56} also follows from
the limit \eqref{8} (yielding Macdonald polynomials as highest degree
part of Koornwinder polynomials) together with the known fact (although
not in the literature) that the
coefficients $u_{\la,\mu}$ ($|\mu|<|\la|$) in \eqref{52} are bounded
as $a_1\to\iy$.
\subsection{Binomial formula for Jack polynomials}
The binomial formula for Jack polynomials was given by
Olounkov \& Olshanski \cite[p.72]{14}:
\begin{equation}
\frac{P_\la(1+x,\tau)}{P_\la(1;\tau)}=
\sum_{\mu\subseteq\la}
\frac{P_\mu^{\rm ip}(\la+\tau\de;\tau)}{P_\mu^{\rm ip}(\mu+\tau\de;\tau)}\,
\frac{P_\mu(x;\tau)}{P_\mu(1;\tau)}\,.
\label{61}
\end{equation}
It is a limit case of \eqref{54} because of the limits
\eqref{22}, \eqref{65} (twice), \eqref{25} (twice) and \eqref{39}.
If we compare \eqref{61} with Macdonald \cite[(6.15), (6.24)]{3},
Lassalle \cite[\S3]{25} or Yan \cite[(10)]{26}, we see that
$P_\mu^{\rm ip}(\la+\tau\de;\tau)/P_\mu^{\rm ip}(\mu+\tau\de;\tau)$
equals the generalized binomial coefficient $\binom\la\mu$ defined in
these references.
\section{Binomial formula for $BC_n$-type Jacobi polynomials}
Let $\al'$ be given by \eqref{81}.
\begin{theorem}
\label{94}
For $BC_n$-type Jacobi polynomials we have the binomial formula
\begin{equation}
\frac{P_\la(x;\tau;\al,\be)}{P_\la(0;\tau;\al,\be)}
=\sum_{\mu\subseteq\la}
\frac{P_\mu^{\rm ip}(\la+\tau\de+\al';\tau,\al')}
{P_\mu^{\rm ip}(\mu+\tau\de+\al';\tau,\al')}\,
\frac{P_\mu(x;\tau)}{P_\mu(0;\tau;\al,\be)}\,.
\label{57}
\end{equation}
\end{theorem}
\Proof
From \eqref{43} we obtain:
\begin{multline*}
\frac{P_\la(x;q,q^\tau;q^{\al+1},-q^{\be+1},1,-1)}
{P_\la(q^{\tau\de+\al+1};q,q^\tau;q^{\al+1},-q^{\be+1},1,-1)}\\
=\sum_{\mu\subseteq\la}
\frac{P_\mu^{\rm ip}(q^{\la+\tau\de+\al'};q,q^\tau,q^{\al'})}
{P_\mu^{\rm ip}(q^{\mu+\tau\de+\al'};q,q^\tau,q^{\al'})}\,
\frac{P_\mu^{\rm ip}(x;q,q^\tau,q^{\al+1})}
{P_\mu(q^{\tau\de+\al+1};q,q^\tau;q^{\al+1},-q^{\be+1},1,-1)}\,.
\end{multline*}
Now let $q\uparrow1$ and apply \eqref{58}, \eqref{60} (twice),
\eqref{31} (twice)
and \eqref{36}.
\qed
\bPP
By comparing the binomial formulas \eqref{57} and \eqref{61} we arrive
at the limit
\begin{equation}
\lim_{\al\to\iy}P_\la(x;\tau;\al,\be)=P_\la(x-1^n;\tau),
\label{62}
\end{equation}
where $1^n$ is the $n$-vector with all coordinates equal to 1.
The limit \eqref{62}
was given slightly more generally in \cite[Theorem 4.2]{10} and goes back to an
unpublished result by Beerends and the author.
For the proof of \eqref{62} let $\al\to\iy$ in the \RHS\ of \eqref{57}
and use \eqref{41}.  Then we will obtain the \RHS\ of \eqref{61} if we
can prove that \eqref{62} is valid for $x=0$.
But in that case \eqref{62} follows by comparing \eqref{59} and \eqref{47}.
\begin{remark}
In \eqref{57} we have an expansion
\begin{equation}
P_\la(x;\tau;\al,\be)=\sum_{\mu\subseteq\la}b_{\la,\mu}\,P_\mu(x;\tau)
\label{64}
\end{equation}
with
\begin{equation}
b_{\la,\mu}=\frac{P_\la(0;\tau;\al,\be)}{P_\mu(0;\tau;\al,\be)}\,
\frac{P_\mu^{\rm ip}(\la+\tau\de+\al';\tau,\al')}
{P_\mu^{\rm ip}(\mu+\tau\de+\al';\tau,\al')}\,,
\label{67}
\end{equation}
where, by \eqref{32},
\begin{equation}
P_\mu^{\rm ip}(\la+\tau\de+\al';\tau,\al')
=\sum_T \psi_T(\tau)\,\prod_{s\in\mu}\big((\la_{T(s)}+\tau\de_{T(s)}+\al')^2
-(\tau\de_{T(s)}+\al'+a_\la'(s)-\tau l_\la'(s))^2\big).
\label{71}
\end{equation}
On the other hand, Macdonald \cite[p.58]{3}
(also in Beerends \& Opdam \cite[(5.12)]{27}) has \eqref{64} with
\begin{equation}
b_{\la,\mu}=(-1)^{|\la|-|\mu|}\frac{P_\la(1;\tau)}{P_\mu(1;\tau)}\,
c_{\la/\mu}(2\tau(n-1)+\al'+1)\,
\prod_{j=1}^n(\mu_j+\tau(n+j-2)+\al+1)_{\la_j-\mu_j}\,.
\label{68}
\end{equation}
Here
\begin{equation}
c_{\la/\mu}(C)=\sum_T f_T(C)
\label{66}
\end{equation}
for a certain function $f_T$ and the sum is over all standard tableaux
$T$ of shape $\la/\mu$. Macdonald derives \eqref{66} by solving a
recurrence relation \cite[(9.16)]{3} for $c_{\la/\mu}(C)$ ($C$ and $\la$
fixed). Thus both in Macdonald's formula \eqref{68} and in formula \eqref{67}
the expansion coefficients $b_{\la,\mu}$ are essentially given
combinatorially by tableau sums (\eqref{66} and \eqref{71}, respectively).
However the tableau sums are quite different: over standard tableaux of
shape $\la/\mu$ in Macdonald's case and over reverse tableaux $T$ of
shape $\mu$ with entries in $\{1,\ldots,n\}$ in the case of
\eqref{71}.
I do not see how to match these tableau sums with each other.
\end{remark}
\section{Explicit expressions for $n=2$}
The combinatorial formulas \eqref{17}, \eqref{18}, \eqref{19} for
Macdonald polynomials
and their ($BC_n$-type) interpolation versions all have the form
\begin{equation}
\sum_T \psi_T(q,t)\,\prod_{s\in\la} f(x_{T(s)}).
\label{83}
\end{equation}
Here $\psi_T(q,t)$ is given by \eqref{24}, while $f$ is an elementary (Laurent)
polynomial which may also depend on $s,T(s),t$ and $a$.
The sum is over all reverse tableaux $T$ of shape $\la$ with entries in
$\{1,\ldots,n\}$.

Now let $n=2$ and $\la:=(m,0)$. Then the possible reverse tableaux $T$
of shape $(m,0)$ with entries in $\{1,2\}$ are the tableaux $T_k$
($k=0,1,\ldots,m$) given by
\begin{equation}
T_k(1,j)=2\quad\mbox{if $j=1,\ldots,k$\quad and\quad
$=1$\quad if $j=k+1,\ldots,m$.}
\label{84}
\end{equation}
So we have to compute $\psi_{T_k}(q,t)$, which is given by \eqref{24} as a
double product involving
$(R\backslash C)_{\la^{(0)}/\la^{(1)}}=\{(1,1),\ldots,(1,k)\}$ and
$(R\backslash C)_{\la^{(1)}/\la^{(2)}}=\emptyset$. Hence
\begin{equation*}
\psi_{T_k}(q,t)=\prod_{s\in (R\backslash C)_{\la^{(0)}/\la^{(1)}}}
\frac{b_{(k)}(s;q,t)}{b_{(m)}(s;q,t)}\,,
\end{equation*}
where $b_\mu(s;q,t)$ is defined by \eqref{82}. This yields for
$j=1,\ldots,k$ that
\begin{equation*}
b_{(k)}((1,j);q,t)=\frac{1-q^{k-j}t}{1-q^{k-j+1}}\,,\qquad
b_{(m)}((1,j);q,t)=\frac{1-q^{m-j}t}{1-q^{m-j+1}}\,.
\end{equation*}
Hence
\begin{equation}
\psi_{T_k}(q,t)=\frac{(q^{k-1}t;q^{-1})_k\,(q^m;q^{-1})_k}
{(q^k;q^{-1})_k\,(q^{m-1}t;q^{-1})_k}
=\frac{(t,q^{-m};q)_k}{(q,q^{1-m}t^{-1};q)_k}\,(qt^{-1})^k.
\label{85}
\end{equation}
\subsection{Explicit expression for
$BC_2$-type interpolation Macdonald polynomials}
Consider in \eqref{83} $\prod_{j=1}^m f(x_{T_k(1,j)})$ for the case of
\eqref{17}. By taking into account \eqref{84} this product becomes
\begin{multline}
\prod_{j=1}^k\big((x_2-q^{j-1}a)(1-(q^{j-1}a)^{-1}x_2^{-1})\big)\,
\prod_{j=k+1}^m\big((x_1-q^{j-1}ta)(1-(q^{j-1}ta)^{-1}x_1^{-1})\big)\\
=(-1)^mq^{-\half m(m-1)}(ta)^{-m}(tax_1,tax_1^{-1};q)_m\,
\frac{(ax_2,ax_2^{-1};q)_k}{(tax_1,tax_1^{-1};q)_k}\,t^k.
\label{86}
\end{multline}
Then take \eqref{83}, a single sum over $T=T_k$ ($k=0,\ldots,m$),
with \eqref{85} and \eqref{86} substituted. There results
\begin{equation}
P_{m,0}^{\rm ip}(x_1,x_2;q,t,a)
=\frac{(tax_1,tax_1^{-1};q)_m}{q^{\half m(m-1)}(-ta)^m}\,
\qhyp43{q^{-m},t,ax_2,ax_2^{-1}}{q^{1-m}t^{-1},tax_1,tax_1^{-1}}{q,q}.
\label{87}
\end{equation}

Next we want to obtain an explicit expression for $P_{m_1,m_2}^{\rm ip}$ by
using the reduction formula \eqref{69}.
This takes in the present case the form
\begin{equation*}
P_{m_1,m_2}^{\rm ip}(x_1,x_2;q,t,a)
=a^{-2m_2} q^{-m_2(m_2-1)}
(ax_1,ax_1^{-1},ax_2,ax_2^{-1};q)_{m_2}\,
P_{m_1-m_2,0}^{\rm ip}(x_1,x_2;q,t,q^{m_2}a).
\end{equation*}
In combination with \eqref{87} this gives
\begin{multline}
P_{m_1,m_2}^{\rm ip}(x_1,x_2;q,t,a)
=\frac{q^{-\half m_1(m_1-1)-\half m_2(m_2-1)}}
{(-t)^{m_1-m_2} a^{m_1+m_2}}\,
(ax_1,ax_1^{-1},ax_2,ax_2^{-1};q)_{m_2}\\
\times(q^{m_2}tax_1,q^{m_2}tax_1^{-1};q)_{m_1-m_2}\,
\qhyp43{q^{-m_1+m_2},t,q^{m_2}ax_2,q^{m_2}ax_2^{-1}}
{q^{1-m_1+m_2}t^{-1},q^{m_2}tax_1,q^{m_2}tax_1^{-1}}{q,q}.
\label{88}
\end{multline}

In order to make clear the symmetry in $x_1,x_2$ in \eqref{88}
we can expand the second line of this formula as
\begin{equation*}
\sum_{k=0}^{m_1-m_2}
\frac{(q^{-m_1+m_2},t;q)_k\,q^k}{(q^{1-m_1+m_2}t^{-1},q:q)_k}\,
(q^{m_2+k}tax_1,q^{m_2+k}tax_1^{-1};q)_{m_1-m_2-k}\,
(q^{m_2}ax_2,q^{m_2}ax_2^{-1};q)_k.
\end{equation*}
Now substitute there a version of the $q$-Pfaff-Saalsch\"utz formula
\cite[(II.12)]{29}:
\begin{equation*}
\frac{(q^{m_2+k}tax_1,q^{m_2+k}tax_1^{-1};q)_{m_1-m_2-k}}
{(q^k t,q^{k+2m_2}ta^2;q)_{m_1-m_2-k}}=
\qhyp32{q^{-m_1+m_2+k},q^{m_2}ax_1,q^{m_2}ax_1^{-1}}
{q^{1-m_1+m_2}t^{-1},q^{2m_2+k}ta^2}{q,q}
\end{equation*}
and expand the ${}_3\phi_2$. Then \eqref{88} takes the form
\begin{multline}
P_{m_1,m_2}^{\rm ip}(x_1,x_2;q,t,a)
=\frac{q^{-\half m_1(m_1-1)-\half m_2(m_2-1)}}
{(-t)^{m_1-m_2} a^{m_1+m_2}}\,
(t,q^{2m_2}ta^2;q)_{m_1-m_2}\\
\times(ax_1,ax_1^{-1},ax_2,ax_2^{-1};q)_{m_2}\,
\sum_{\substack{j,k\ge0\\j+k\le m_1-m_2}}
\frac{(q^{-m_1+m_2};q)_{j+k}q^{j+k}}{(q^{2m_2}ta^2;q)_{j+k}}
\\
\times
\frac{(q^{m_2}ax_1,q^{m_2}ax_1^{-1};q)_j}{(q^{1-m_1+m_2}t^{-1},q;q)_j}\,
\frac{(q^{m_2}ax_2,q^{m_2}ax_2^{-1};q)_k}{(q^{1-m_1+m_2}t^{-1},q;q)_k}\,.
\end{multline}
\subsection{(Interpolation) Macdonald polynomials
for~\mbox{$n=2$}}
For explicit formulas of $P_{m_1,m_2}^{\rm ip}(x_1,x_2;q,t)$ and
$P_{m_1,m_2}(x_1,x_2;q,t)$ we may either give a derivation analogous
to the one for \eqref{88}, now starting from \eqref{18} or \eqref{17}, or apply,
much quicker, the limits \eqref{21} or \eqref{20} to \eqref{88}.
For $P_{m_1,m_2}^{\rm ip}(x_1,x_2;q,t)$ we obtain
\begin{multline*}
P_{m_1,m_2}^{\rm ip}(x_1,x_2;q,t)=x_1^{m_1} x_2^{m_2}\,
(x_1^{-1},x_2^{-1};q)_{m_2} (q^{m_2}tx_1^{-1};q)_{m_1-m_2}\\
\times \qhyp32{q^{-m_1+m_2},t,q^{m_2}x_2^{-1}}
{q^{1-m_1+m_2}t^{-1},q^{m_2}tx_1^{-1}}{q,\frac{qx_2}{tx_1}}.
\end{multline*}
By series inversion this yields:
\begin{equation}
P_{m_1,m_2}^{\rm ip}(x_1,x_2;q,t)=x_1^{m_2} x_2^{m_1}\,(x_1^{-1};q)_{m_2}\,
(x_2^{-1};q)_{m_1}\,
\qhyp32{q^{-m_1+m_2},t,q^{-m_1+1}t^{-1}x_1}{q^{1-m_1+m_2}t^{-1},q^{-m_1+1}x_2}
{q,q}.
\label{89}
\end{equation}
Different explicit expansions for interpolation Macdonald polynomials in case
$n=2$ are given by
Morse \cite[Theorem 1.1]{34}, \cite[Theorem 1]{33}.

For $P_{m_1,m_2}(x_1,x_2;q,t)$ we obtain:
\begin{align}
P_{m_1,m_2}(x_1,x_2;q,t)&=x_1^{m_1} x_2^{m_2}\,
\qhyp21{q^{-m_1+m_2},t}{q^{1-m_1+m_2}t^{-1}}{q,\frac{qx_2}{tx_1}}
\label{90}\\
&=\frac{(q;q)_{m_1-m_2}}{(t;q)_{m_1-m_2}}\,
\sum_{j=0}^{m_1-m_2}\frac{(t;q)_j(t;q)_{m_1-m_2-j}}{(q;q)_j(q;q)_{m_1-m_2-j}}\,
x_1^{m_1-j}x_2^{m_2+j}
\label{91}\\
&=\frac{(q;q)_{m_1-m_2}}{(t;q)_{m_1-m_2}}\,(x_1x_2)^{\half(m_1+m_2)}\,
C_{m_1-m_2}\left(\frac{x_1+x_2}{2(x_1x_2)^\half};t\mid q\right).
\label{92}
\end{align}
Here we used the {\em $q$-ultraspherical polynomial}
(see \cite[(4.4)]{31}):
\begin{equation}
C_m(\cos\tha;t\mid q):=
\sum_{j=0}^m\frac{(t;q)_j(t;q)_{m-j}}{(q;q)_j(q;q)_{m-j}}\,
e^{i(m-2j)\tha}.
\end{equation}
Formula \eqref{92} was earlier given by Morse \cite[Theorem 2]{33}.
\subsection{Jack polynomials and ($BC_n$) interpolation Jack polynomials
for~\mbox{$n=2$}}
We can get explicit formulas for the $n=2$ cases of $BC_n$-type interpolation
Jack polynomials, interpolation Jack polynomials and Jack polynomials by
applying the limit formulas \eqref{31}, \eqref{25}, \eqref{22}
to \eqref{88}, \eqref{89}, \eqref{90}, respectively. We obtain:
\begin{multline}
P_{m_1,m_2}^{\rm ip}(x_1,x_2;\tau,\al)=(-1)^{m_1+m_2}\,
(\al+x_1,\al-x_1,\al+x_2,\al-x_2)_{m_2}\\[\smallskipamount]
\times(m_2+\tau+\al+x_1,m_2+\tau+\al-x_1)_{m_1-m_2}\\[\smallskipamount]
\times \hyp43{-m_1+m_2,\tau,m_2+\al+x_2,m_2+\al-x_2}
{1-m_1+m_2-\tau,m_2+\tau+\al+x_1,m_2+\tau+\al-x_1}1,
\end{multline}
\begin{equation}
P_{m_1,m_2}^{\rm ip}(x_1,x_2;\tau)=(-1)^{m_1+m_2}(-x_1)_{m_2} (-x_2)_{m_1}\,
\hyp32{-m_1+m_2,\tau,-m_1+1-\tau+x_1}{1-m_1+m_2-\tau,-m_1+1+x_2}1.
\end{equation}
\begin{align}
P_{m_1,m_2}(x_1,x_2;\tau)
&=x_1^{m_1}x_2^{m_2}\,\hyp21{-m_1+m_2,\tau}{1-m_1+m_2-\tau}{\frac{x_2}{x_1}}
\nonumber\\
&=\frac{(m_1-m_2)!}{(\tau)_{m_1-m_2}}\,
\sum_{j=0}^{m_1-m_2}\frac{(\tau)_j(\tau)_{m_1-m_2-j}}{j!\,(m_1-m_2-j)!}\,
x_1^{m_1-j}x_2^{m_2+j}\nonumber\\
&=\frac{(m_1-m_2)!}{(\tau)_{m_1-m_2}}\,(x_1x_2)^{\half(m_1+m_2)}\,
C_{m_1-m_2}^\tau\left(\frac{x_1+x_2}{2(x_1x_2)^\half}\right).
\label{96}
\end{align}
Here we used the {\em ultraspherical polynomial}
(see \cite[Section 10.9]{32}):
\begin{equation}
C_m^\tau(\cos\tha):=
\sum_{j=0}^m\frac{(\tau)_j(\tau)_{m-j}}{j!\,(m-j)!}\,e^{i(m-2j)\tha}.
\end{equation}
\section{Binomial formulas for $n=2$}
In this section we will only discuss the $n=2$ case of the binomial
formulas for $BC$-type polynomials (Koornwinder and Jacobi), because
they lead to explicit expressions of these polynomials. We will not
discuss here the $n=2$ cases of the binomial formulas for the
Macdonald and Jack polynomials, although they have interesting aspects
from the point of view of special functions.
\subsection{Binomial formula for Koornwinder polynomials for $n=2$}
The binomial formula \eqref{43} takes for $n=2$ the form
\begin{multline} 
\frac{P_{m_1,m_2}(x_1,x_2;q,t;a_1,a_2,a_3,a_4)}
{P_{m_1,m_2}(ta_1,a_1;q,t;a_1,a_2,a_3,a_4)}
=\sum_{k_2=0}^{m_2}\sum_{k_1=k_2}^{m_1}
\frac{P_{k_1,k_2}^{\rm ip}(q^{m_1}ta_1',q^{m_2}a_1';q,t,a_1')}
{P_{k_1,k_2}^{\rm ip}(q^{k_1}ta_1',q^{k_2}a_1';q,t,a_1')}\\
\times\frac{P_{k_1,k_2}^{\rm ip}(x_1,x_2;q,t,a_1)}
{P_{k_1,k_2}(ta_1,a_1;q,t;a_1,a_2,a_3,a_4)}\,.
\label{93}
\end{multline}
This gives an explicit expression for Koornwinder polynomials for $n=2$,
since we have explicit expressions for everything on the \RHS.
For the first quotient on the \RHS\ of \eqref{93} we get by \eqref{88}:
\begin{multline}
\frac{P_{k_1,k_2}^{\rm ip}(q^{m_1}ta_1',q^{m_2}a_1';q,t,a_1')}
{P_{k_1,k_2}^{\rm ip}(q^{k_1}ta_1',q^{k_2}a_1';q,t,a_1')}
=\frac{(q^{m_1}ta_1'^2,q^{-m_1}t^{-1},q^{m_2}a_1'^2,q^{-m_2};q)_{k_2}}
{(q^{k_1}ta_1'^2,q^{-k_1}t^{-1},q^{k_2}a_1'^2,q^{-k_2};q)_{k_2}}\\
\times\frac{(q^{m_1+k_2}t^2a_1'^2,q^{m_2-k_1};q)_{k_1-k_2}}
{(q^{k_1+k_2}t^2a_1'^2,q^{k_2-k_1};q)_{k_1-k_2}}\,
\qhyp43{q^{-k_1+k _2},t,q^{m_2+k_2}a_1'^2,q^{m_2-k_2}}
{q^{1-k_1+k_2}t^{-1},q^{m_2+k_2}t^2a_1'^2,q^{-m_1+k_2}}{q,q}.
\label{95}
\end{multline}
As for the second quotient on the \RHS\  the numerator is obtained
from \eqref{88}, i.e.,
\begin{multline}
P_{k_1,k_2}^{\rm ip}(x_1,x_2;q,t,a_1)
=\frac{q^{-\half k_1(k_1-1)-\half k_2(k_2-1)}}
{(-t)^{k_1-k_2} a_1^{k_1+k_2}}\,
(a_1x_1,a_1x_1^{-1},a_1x_2,a_2x_2^{-1};q)_{k_2}\\
\times(q^{k_2}ta_1x_1,q^{k_2}ta_1x_1^{-1};q)_{k_1-k_2}\,
\qhyp43{q^{-k_1+k_2},t,q^{k_2}a_1x_2,q^{k_2}a_1x_2^{-1}}
{q^{1-k_1+k_2}t^{-1},q^{k_2}ta_1x_1,q^{k_2}ta_1x_1^{-1}}{q,q},
\end{multline}
and the denominator is obtained from \eqref{48}, i.e.,
\begin{multline}
P_{k_1,k_2}(ta_1,a_1;q,t;a_1,a_2,a_3,a_4)=
t^{-k_1}a_1^{-k_1+k_2}\,
\frac{(ta_1'^2;q)_{k_1}\,(a_1'^2;q)_{k_2}}
{(t^2a_1'^2;q)_{2k_1}\,(a_1'^2;q)_{2k_2}}\\
\times(ta_1a_2,ta_1a_3,ta_1a_4;q)_{k_1}\,
(a_1a_2,a_1a_3,a_1a_4;q)_{k_2}\,
\frac{(t^2a_1'^2;q)_{k_1+k_2}}{(ta_1'^2;q)_{k_1+k_2}}\,
\frac{(t^2;q)_{k_1-k_2}}{(t;q)_{k_1-k_2}}\,.
\end{multline}
\subsection{Binomial formula for $BC_2$-type Jacobi polynomials}
We can specialize the binomial formula \eqref{57} for $BC_n$-type
Jacobi polynomials to $n=2$ and get everything on the \RHS\ explicitly:
\begin{equation}
\frac{P_{m_1,m_2}(x_1,x_2;\tau;\al,\be)}{P_{m_1,m_2}(0,0;\tau;\al,\be)}
=\sum_{k_2=0}^{m_2}\sum_{k_1=k_2}^{m_1}
\frac{P_{k_1,k_2}^{\rm ip}(m_1+\tau+\al',m_2+\al';\tau,\al')}
{P_{k_1,k_2}^{\rm ip}(k_1+\tau+\al',k_2+\al';\tau,\al')}\,
\frac{P_{k_1,k_2}(x_1,x_2;\tau)}{P_{k_1,k_2}(0,0;\tau;\al,\be)}\,.
\label{97}
\end{equation}
The first quotient on the \RHS\ can be explicitly evaluated as a limit case of
\eqref{95} by using \eqref{31}:
\begin{multline}
\frac{P_{k_1,k_2}^{\rm ip}(m_1+\tau+\al',m_2+\al';\tau,\al')}
{P_{k_1,k_2}^{\rm ip}(k_1+\tau+\al',k_2+\al';\tau,\al')}=
\frac{(m_1+\tau+2\al',-m_1-\tau,m_2+2\al',-m_2)_{k_2}}
{(k_1+\tau+2\al',-k_1-\tau,k_2+2\al',-k_2)_{k_2}}\\
\times\frac{(m_2+k_2+2\tau+2\al',m_2-k_2)_{k_1-k_2}}
{(k_1+k_2+2\tau+2\al',k_2-k_1)_{k_1-k_2}}\\
\times\hyp43{-k_1+k_2,\tau,m_2+k_2+2\al',m_2-k_2}
{1-k_1+k_2-\tau,m_2+k_2+2\tau+2\al',-m_1+k_2}1.
\label{98}
\end{multline}
For the numerator of the second quotient we have \eqref{96}:
\begin{equation}
P_{k_1,k_2}(x_1,x_2;\tau)=
\frac{(q;q)_{k_1-k_2}}{(t;q)_{k_1-k_2}}\,(x_1x_2)^{\half(k_1+k_2)}\,
C_{k_1-k_2}^\tau\left(\frac{x_1+x_2}{2(x_1x_2)^\half}\right).
\label{99}
\end{equation}
The denominator of the second quotient can be obtained by specializing
\eqref{63}:
\begin{multline}
P_{k_1,k_2}(0,0;\tau;\al,\be)=(-1)^{k_1+k_2}\,
\frac{(\tau+2\al',\tau+\al+1)_{k_1}\,(2\al',\al+1)_{k_2}}
{(2\tau+2\al')_{2k_1}\,(2\al')_{2k_2}}\\
\times\frac{(2\tau+2\al')_{k_1+k_2}}{(\tau+2\al')_{k_1+k_2}}\,
\frac{(2\tau)_{k_1-k_2}}{(\tau)_{k_1-k_2}}\,.
\label{100}
\end{multline}
The explicit formula for $BC_2$-type Jacobi polynomials which is given by
combination of \eqref{97}, \eqref{98}, \eqref{99} and \eqref{100} was obtained
in a very different way by Koornwinder \& Sprinkhuizen
\cite[Corollary 6.6]{4}.
\quad\\
\begin{footnotesize}
\begin{quote}
T. H. Koornwinder, Korteweg-de Vries Institute, University of
 Amsterdam,\\
 P.O.\ Box 94248, 1090 GE Amsterdam, The Netherlands;
 
\vspace{\smallskipamount}
email: {\tt T.H.Koornwinder@uva.nl}
\end{quote}
\end{footnotesize}

\end{document}